\documentclass{birkjour}

\usepackage{amsmath,amssymb,amsxtra,latexsym,amsthm,amscd,amsfonts,mathrsfs,pb-diagram}
\usepackage[colorlinks=true,linkcolor=magenta,citecolor=blue]{hyperref}
\usepackage{mathabx,bezier}
\usepackage{enumitem}
\usepackage[latin1]{inputenc}

\newtheorem{dn}{Definition}[section]
\newtheorem{dl}{Theorem}[section]
\newtheorem{md}{Proposition}[section]
\newtheorem{bd}{Lemma}[section]
\newtheorem{hq}{Corollary}[section]
\newtheorem{nx}{Remark}[section]
\newtheorem{vd}{Example}[section]
\newcommand{\R}{\mathbb{R}}
\newcommand{\C}{\mathbb{C}}
\newcommand{\Z}{\mathbb{Z}}

\newcommand{\ity}{\infty}
\newcommand{\f}{\frac}

\newcommand{\bbd}{\begin{bd}}
\newcommand{\ebd}{\end{bd}}
\newcommand{\bdn}{\begin{dn}}
\newcommand{\edn}{\end{dn}}
\newcommand{\bhq}{\begin{hq}}
\newcommand{\ehq}{\end{hq}}
\newcommand{\bdl}{\begin{dl}}
\newcommand{\edl}{\end{dl}}
\newcommand{\bnx}{\begin{nx}}
\newcommand{\enx}{\end{nx}}
\newcommand{\bmd}{\begin{md}}
\newcommand{\emd}{\end{md}}
\newcommand{\bvd}{\begin{vd}}
\newcommand{\evd}{\end{vd}}

\begin{document}

\title[\small{$L^1$-$L^1$ estimates for solutions to visco-elastic damped $\sigma$-evolution models}]{Some $L^1$-$L^1$ estimates for solutions to visco-elastic damped $\sigma$-evolution models}
\author{Tuan Anh Dao}
\address{
$\quad$ School of Applied Mathematics and Informatics, \\
Hanoi University of Science and Technology, \\
No.1 Dai Co Viet road, Hanoi, Vietnam.
\hfill\break
Faculty for Mathematics and Computer Science, \\
TU Bergakademie Freiberg, \\
Pr\"{u}ferstr. 9, 09596, Freiberg, Germany.}
\email{anh.daotuan@hust.edu.vn}

\begin{abstract}
This note is to conclude $L^1-L^1$ estimates for solutions to the following Cauchy problem for visco-elastic damped $\sigma$-evolution models:
\begin{equation}
\begin{cases}
u_{tt}+ (-\Delta)^\sigma u+ (-\Delta)^\sigma u_t = 0, &\quad x\in \R^n,\, t \ge 0, \\
u(0,x)= u_0(x),\quad u_t(0,x)=u_1(x), &\quad x\in \R^n, \label{pt1.1}
\end{cases}
\end{equation}
where $\sigma> 1$, in all space dimensions $n\ge 1$.
\end{abstract}
\maketitle

\section{Introduction}
At prsent, there has been very little work on the question of getting $L^1-L^1$ estimates for solutions to (\ref{pt1.1}). Back in $2000$, let us first recall the pioneering paper of Shibata \cite{Shibata} devoting to the study of one of the most well-known equations of (\ref{pt1.1}) in the case $\sigma=1$, the so-called strongly damped wave equation. In the cited paper, relying on the very special structure of fundamental solutions to the wave equation he succeeded in obtaining the following $L^1-L^1$ estimates:
$$ \|u(t,\cdot)\|_{L^1}\lesssim
\begin{cases}
(1+t)^{\frac{n}{4}}\|u_0\|_{L^1}+ (1+t)^{\frac{n+2}{4}}\|u_1\|_{L^1} &\text{ if }n\ge 2 \text{ is even}, \\
(1+t)^{\frac{n-1}{4}}\|u_0\|_{L^1}+ (1+t)^{\frac{n+1}{4}}\|u_1\|_{L^1} &\text{ if }n\ge 3 \text{ is odd}, \\
\end{cases}$$ 
by taking into considerations the connection to Fourier multipliers appearing for wave models. Quite recently, regarding a different interesting model of (\ref{pt1.1}), namely that with $\sigma=2$, D'Abbicco and his collaborators \cite{DabbiccoGirardiLiang} have employed a different technique in comparison with that used in \cite{Shibata} to derive $L^1-L^1$ estimates for solutions to the strongly damped plate equation as follows:
$$ \|u(t,\cdot)\|_{L^1}\lesssim (1+t)^{\frac{n}{4}}\|u_0\|_{L^1}+ (1+t)^{\frac{n+2}{4}}\|u_1\|_{L^1} \quad \text{ for any }n\ge 5. $$
Here this limitation of the space dimensions comes from the technical difficulty. More precisely, the authors used Bernstein inequality to estimate Fourier multipliers for small frequencies. To apply this technique, it is necessary to require the above restriction to space dimensions. Independently from the above mentioned results, Dao-Reissig in the recent paper \cite{DaoReissig} have considered the more general cases of (\ref{pt1.1}) for any $\sigma> 1$ to achieve $L^1- L^1$ estimates for solutions localized to small frequencies by applying the theory of modified Bessel functions linked to Fa\`{a} di Bruno's formula. Unfortunately, this strategy fails in the treatment of large frequencies. For this reason, the present paper is to fill this lack and report some $L^1- L^1$ estimates for solutions to (\ref{pt1.1}) as well.
\par In order to state our main result, we introduce the following notations used in this paper:
\begin{itemize}[leftmargin=*]
\item We write $f\lesssim g$ when there exists a constant $C>0$ such that $f\le Cg$, and $f \approx g$ when $g\lesssim f\lesssim g$.
\item The spaces $H^a_1$ with $a \ge 0$ stand for Bessel potential spaces based on $L^1$ spaces, where $\big<D\big>^a$ denote the pseudo-differential operators with symbols $\big<\xi\big>^a$. 
\item For a given number $s \in \R$, we denote $[s]:= \max \big\{k \in \Z \,\, : \,\, k\le s \big\}$ and $[s]^+:= \max\{s,0\}$ as its integer part and its positive part, respectively.
\end{itemize}
\par The main purpose of this note is to prove the following result.
\bdl[Main result] \label{dl1.1}
Let $\sigma>1$. Then, the Sobolev solutions to (\ref{pt1.1}) satisfy the following $L^1- L^1$ estimates:
\begin{align*}
\big\||D|^a u (t,\cdot)\big\|_{L^1}& \lesssim (1+t)^{\frac{1}{2}(2+[\frac{n}{2}])-\frac{a}{2\sigma}}\|u_0\|_{H^a_1} \\
&\quad + (1+t)^{1+\frac{1}{2}(1+[\frac{n}{2}])-\frac{a}{2\sigma}}\|u_1\|_{H^{[a-\sigma]^+}_1}, \\
\big\||D|^a u_t (t,\cdot)\big\|_{L^1}& \lesssim (1+t)^{\frac{1}{2}(1+[\frac{n}{2}])-\frac{a}{2\sigma}}\|u_0\|_{H^{2\sigma+ [a-\sigma]^+}_1} \\
&\quad + (1+t)^{\frac{1}{2}(2+[\frac{n}{2}])-\frac{a}{2\sigma}}\|u_1\|_{H^{2\sigma+ [a-\sigma]^+}_1},
\end{align*}
where $a \ge 0$ and for all space dimensions $n\ge 1$.
\edl

\bnx
\fontshape{n}
\selectfont
Here we want to underline that at the first glance the decay estimates for solutions produced from the results of \cite{Shibata} or \cite{DabbiccoGirardiLiang} are somehow better than those of Theorem \ref{dl1.1} when we choose $\sigma=1$ formally or $\sigma=2$, respectively. The fact is that this comes from the very special structure of solutions to (\ref{pt1.1}) in the cases of $\sigma=1$ and $\sigma=2$ coupled with some used techniques under the suitable restriction to space dimensions. However, these obtained results from Theorem \ref{dl1.1} imply an important regconition that we may conclude the desired $L^1- L^1$ estimates for solutions to (\ref{pt1.1}) in gerenal cases of $\sigma>1$ without any constraint condition to space dimensions, i.e. our main result is valid for any $n\ge 1$.
\enx
	
\section{Proof of main result} \label{Sec.Prof}
At first, using partial Fourier transformation to (\ref{pt1.1}) we obtain the Cauchy problem for $\widehat{u}(t,\xi):= \mathfrak{F}(u(t,x)),$ $\widehat{u}_0(\xi):= \mathfrak{F}(u_0(x))$ and $\widehat{u}_1(\xi):= \mathfrak{F}(u_1(x))$ as follows:
\begin{equation}
\widehat{u}_{tt}+ |\xi|^{2\sigma} \widehat{u}_t+ |\xi|^{2\sigma} \widehat{u}=0,\quad \widehat{u}(0,\xi)= \widehat{u}_0(\xi),\quad \widehat{u}_t(0,\xi)=\widehat{u}_1(\xi).
\label{pt1.2}
\end{equation}
The characteristic roots are
$$ \lambda_{1,2}=\lambda_{1,2}(\xi)= \f{1}{2}\Big(-|\xi|^{2\sigma}\pm \sqrt{|\xi|^{4\sigma}-4|\xi|^{2\sigma}}\Big). $$
The solutions to (\ref{pt1.2}) are presented by the following formula (here we assume $\lambda_1 \neq \lambda_2$):
\begin{align*}
\widehat{u}(t,\xi)&= \frac{\lambda_1 e^{\lambda_2 t}-\lambda_2 e^{\lambda_1 t}}{\lambda_1- \lambda_2}\widehat{u}_0(\xi)+ \frac{e^{\lambda_1 t}-e^{\lambda_2 t}}{\lambda_1- \lambda_2}\widehat{u}_1(\xi) \\ 
&=: \widehat{\mathcal{K}}_0(t,\xi)\widehat{u}_0(\xi)+ \widehat{\mathcal{K}}_1(t,\xi)\widehat{u}_1(\xi).
\end{align*}
Taking account of the cases of small and large frequencies separately we have
\begin{itemize}[leftmargin=*]
\item[1.] $\lambda_{1,2}=\lambda_{1,2}(\xi)= -\f{1}{2}\big(|\xi|^{2\sigma}\mp i\sqrt{4|\xi|^{2\sigma}-|\xi|^{4\sigma}}\big)$ \\
and $\lambda_{1,2}\sim -|\xi|^{2\sigma}\pm i|\xi|^\sigma, \quad \lambda_1-\lambda_2 \sim i|\xi|^\sigma \quad$ for $|\xi| \in (0,4^{-\frac{1}{\sigma}})$,
\item[2.] $\lambda_{1,2}=\lambda_{1,2}(\xi)= -\f{1}{2}\big(|\xi|^{2\sigma}\mp \sqrt{|\xi|^{4\sigma}-4|\xi|^{2\sigma}}\big)$ \\
and $\lambda_1\sim -1, \quad \lambda_2\sim -|\xi|^{2\sigma}, \quad \lambda_1-\lambda_2 \sim |\xi|^{2\sigma} \quad$ for $|\xi| \in (4^{\frac{1}{\sigma}},\ity)$.
\end{itemize}
Let $\chi_k= \chi_k(|\xi|)$ with $k=1,2,3$ be smooth cut-off functions having the following properties:
\begin{align*}
&\chi_1(|\xi|)=
\begin{cases}
1 &\quad \text{ if }|\xi|\le 4^{-\frac{1}{\sigma}}, \\
0 &\quad \text{ if }|\xi|\ge 3^{-\frac{1}{\sigma}},
\end{cases}
\qquad
\chi_3(|\xi|)=
\begin{cases}
1 &\quad \text{ if }|\xi|\ge 4^{\frac{1}{\sigma}}, \\
0 &\quad \text{ if }|\xi|\le 3^{\frac{1}{\sigma}},
\end{cases} \\ 
&\text{and } \chi_2(|\xi|)= 1- \chi_1(|\xi|)- \chi_3(|\xi|).
\end{align*}
We note that $\chi_2(|\xi|)= 1$ if $3^{-\frac{1}{\sigma}}\le |\xi|\le 3^{\frac{1}{\sigma}}$ and $\chi_2(|\xi|)= 0$ if $|\xi| \le 4^{-\frac{1}{\sigma}}$ or $|\xi| \ge 4^{\frac{1}{\sigma}}$. Let us now decompose the solutions to (\ref{pt1.1}) into three parts localized individually to small, middle and large frequencies, that is,
$$ u(t,x)= u_{\chi_1}(t,x)+ u_{\chi_2}(t,x)+ u_{\chi_3}(t,x), $$
where
$$u_{\chi_k}(t,x)= \mathfrak{F}^{-1}\big(\chi_k(|\xi|)\widehat{u}(t,\xi)\big)\quad \text{ with } k=1,2,3. $$
For this reason, we shall divide our considerations into three cases as follows.

\subsection{Estimates for small frequencies}
We follow the staments from Corollary $3.3$ in the paper of Dao-Reissig \cite{DaoReissig} to obtain the following estimates for small frequencies.
\bmd
Let $\sigma>1$ and $n\ge 1$. The Sobolev solutions to (\ref{pt1.1}) satisfy the $L^1- L^1$ estimates
\begin{align*}
\big\||D|^a u_{\chi_1} (t,\cdot)\big\|_{L^1}& \lesssim (1+t)^{\frac{1}{2}(2+[\frac{n}{2}])-\frac{a}{2\sigma}}\|u_0\|_{L^1} \\
&\qquad + (1+t)^{1+\frac{1}{2}(1+[\frac{n}{2}])-\frac{a}{2\sigma}}\|u_1\|_{L^1}, \\
\big\||D|^a \partial_t u_{\chi_1} (t,\cdot)\big\|_{L^1}& \lesssim (1+t)^{\frac{1}{2}(1+[\frac{n}{2}])-\frac{a}{2\sigma}}\|u_0\|_{L^1} \\
&\qquad + (1+t)^{\frac{1}{2}(2+[\frac{n}{2}])-\frac{a}{2\sigma}}\|u_1\|_{L^1},
\end{align*}
for any $a \ge 0$.
\label{md2.1}
\emd

\subsection{Estimates for large frequencies}
At first, let us represent the characteristic roots in the form
\begin{equation}
\lambda_1(\xi)= -1-\phi(\xi) \text{ and } \lambda_2(\xi)= -|\xi|^{2\sigma}+1+\phi(\xi), \label{characteristicRoot1}
\end{equation}
where
\begin{equation}
\phi(\xi)=-1+ \int_0^{1}\Big(1- \frac{4}{|\xi|^{2\sigma}}s\Big)^{-\frac{1}{2}}ds. \label{characteristicRoot2}
\end{equation}
For the sake of transparent representation for large frequencies, we introduce the following notations:
\begin{align*}
&\mathcal{K}^1_{u_0}(t,x):= \mathfrak{F}^{-1}\Big(\frac{\lambda_2(\xi) e^{\lambda_1(\xi) t}}{\lambda_1(\xi)- \lambda_2(\xi)}\widehat{u}_0(\xi)\chi_3(|\xi|) \Big)(t,x), \\
&\mathcal{K}^2_{u_0}(t,x):= \mathfrak{F}^{-1}\Big(\frac{\lambda_1(\xi) e^{\lambda_2(\xi) t}}{\lambda_1(\xi)- \lambda_2(\xi)}\widehat{u}_0(\xi)\chi_3(|\xi|) \Big)(t,x), \\
&\mathcal{K}^1_{u_1}(t,x):= \mathfrak{F}^{-1}\Big(\frac{e^{\lambda_1(\xi) t}}{\lambda_1(\xi)- \lambda_2(\xi)}\widehat{u}_1(\xi)\chi_3(|\xi|) \Big)(t,x), \\
&\mathcal{K}^2_{u_1}(t,x):= \mathfrak{F}^{-1}\Big(\frac{e^{\lambda_2(\xi) t}}{\lambda_1(\xi)- \lambda_2(\xi)}\widehat{u}_1(\xi)\chi_3(|\xi|) \Big)(t,x).
\end{align*}
Then, our main goal of this section is to show the following assertions.

\bmd
Let $\sigma>1$ and $n\ge 1$. The following estimates hold:
\begin{align*}
&\big\|\partial_t^{j}|D|^a \mathcal{K}^1_{u_0}(t,\cdot)\big\|_{L^1} \lesssim e^{-ct} \|u_0\|_{H^a_1}, \\
&\big\|\partial_t^{j}|D|^a \mathcal{K}^2_{u_0}(t,\cdot)\big\|_{L^1} \lesssim e^{-ct} \|u_0\|_{H^{2\sigma j+ [a-\sigma]^+}_1}, \\
&\big\|\partial_t^{j}|D|^a \mathcal{K}^1_{u_1}(t,\cdot)\big\|_{L^1} \lesssim e^{-ct} \|u_1\|_{H^{[a-\sigma]^+}_1}, \\
&\big\|\partial_t^{j}|D|^a \mathcal{K}^2_{u_1}(t,\cdot)\big\|_{L^1} \lesssim e^{-ct} \|u_1\|_{H^{2\sigma j+ [a-\sigma]^+}_1},
\end{align*}
where $c$ is a suitable positive constant, for any $t> 0$, $a \ge 0$ and for all integer number $j \ge 0$.
\label{md2.21}
\emd
\noindent In order to prove Proposition \ref{md2.21} let us recall the following auxiliary estimates from Lemma 3.5 in \cite{DaoReissig}.
\bbd \label{bd0.2.1}
The following estimates hold in $\R^n$ for sufficiently large $|\xi|$:
\begin{align}
&\big|\partial_\xi^\alpha |\xi|^{2p\sigma} \big| \lesssim |\xi|^{2p\sigma-|\alpha|}\text{ for all }\alpha\text{ and } p\in\R, \label{bd0.2.1.11} \\
&\big|\partial_\xi^\alpha \phi(\xi)\big| \lesssim |\xi|^{-2\sigma-|\alpha|}\text{ for all }\alpha, \label{bd0.2.1.3} \\
&\Big|\partial_\xi^\alpha \Big(\frac{\lambda_1(\xi) e^{\lambda_2(\xi)t} \lambda_2^j(\xi) |\xi|^b}{\lambda_1(\xi)- \lambda_2(\xi)} \Big) \Big| \lesssim e^{-ct}|\xi|^{2\sigma j+b-2\sigma-|\alpha|} \text{ for all }\alpha, \label{bd0.2.1.14} \\
&\text{ for any } b\in \R,\, j\ge 0 \text{ and } t>0, \text{ where $c$ is a suitable positive constant}, \nonumber
\end{align}
\begin{align}
&\Big|\partial_\xi^\alpha \Big(\frac{e^{\lambda_2(\xi)t} \lambda_2^j(\xi) |\xi|^b}{\lambda_1(\xi)- \lambda_2(\xi)} \Big) \Big| \lesssim e^{-ct}|\xi|^{2\sigma j+b-2\sigma-|\alpha|} \text{ for all }\alpha, \label{bd0.2.1.15} \\
&\text{ for any } b\in \R,\, j\ge 0 \text{ and } t>0, \text{ where $c$ is a suitable positive constant}, \nonumber \\
&\Big|\partial_\xi^\alpha \Big(\frac{\lambda_2(\xi) e^{\lambda_1(\xi)t} \lambda_1^j(\xi) |\xi|^b }{\lambda_1(\xi)- \lambda_2(\xi)} \Big) \Big| \lesssim e^{-ct}|\xi|^{b-|\alpha|}\text{ for all }\alpha,  \label{bd0.2.1.16} \\
&\text{ for any } b\in \R,\, j\ge 0 \text{ and } t>0, \text{ where $c$ is a suitable positive constant}, \nonumber \\
&\Big|\partial_\xi^\alpha\Big(\frac{e^{\lambda_1(\xi)t} \lambda_1^j(\xi) |\xi|^b}{\lambda_1(\xi)- \lambda_2(\xi)} \Big) \Big| \lesssim e^{-ct}|\xi|^{b-2\sigma-|\alpha|} \text{ for all }\alpha, \label{bd0.2.1.17} \\
&\text{ for any } b\in \R,\, j\ge 0 \text{ and } t>0, \text{ where $c$ is a suitable positive constant}. \nonumber
\end{align}
\ebd

\begin{proof}[Proof of Proposition \ref{md2.21}]
In order to indicate some estimates for $\mathcal{K}^2_{u_0}(t,x)$, we may wirte
\begin{align*}
&\partial_t^{j}|D|^a \mathcal{K}^2_{u_0}(t,x) \\
&\quad= \mathfrak{F}^{-1}\Big(\frac{ \lambda_1(\xi) e^{\lambda_2(\xi)t} \lambda_2^j(\xi) |\xi|^{\min\{a,\sigma\}- 2\sigma j}}{\lambda_1(\xi)- \lambda_2(\xi)}\chi_3(|\xi|) |\xi|^{2\sigma j +[a-\sigma]^+}\widehat{u}_0(\xi) \Big)(t,x) \\
&\quad= \mathfrak{F}^{-1}\Big(\frac{ \lambda_1(\xi) e^{\lambda_2(\xi)t} \lambda_2^j(\xi) |\xi|^{\min\{a,\sigma\}- 2\sigma j}}{\lambda_1(\xi)- \lambda_2(\xi)}\chi_3(|\xi|)\Big)(t,x) \ast  |D|^{2\sigma j +[a-\sigma]^+}u_0(x) \\
&\quad=: \mathfrak{F}^{-1}\big(\widehat{\mathcal{L}}^2_0(t,\xi)\big)(t,x) \ast  |D|^{2\sigma j +[a-\sigma]^+}u_0(x).
\end{align*}
By choosing $b= \min\{a,\sigma\}- 2\sigma j$ in (\ref{bd0.2.1.14}), we get
$$\Big|\partial_\xi^\alpha \big(\widehat{\mathcal{L}}^2_0(t,\xi) \big)\Big| \lesssim e^{-ct}|\xi|^{\min\{a,\sigma\}- 2\sigma- |\alpha|} \lesssim e^{-ct}|\xi|^{-\sigma- |\alpha|}, $$
where $c$ is a suitable positive constant. Since
\begin{equation}
e^{ix\xi}= \sum_{k=1}^n \frac{x_k}{i|x|^2} \partial_{\xi_k} e^{ix\xi}, \label{Fouriertransform}
\end{equation}
carrying out $m$ steps of partial integration we derive
$$ \mathfrak{F}^{-1}\big(\widehat{\mathcal{L}}^2_0(t,\xi)\big)(t,x)= C \sum_{|\alpha|=m} \frac{(ix)^\alpha}{|x|^{2|\alpha|}} \mathfrak{F}^{-1}\Big(\partial_\xi^\alpha \big(\widehat{\mathcal{L}}^2_0(t,\xi) \big) \Big)(t,x). $$
For this reason, we obtain the following estimates:
\begin{align*}
\big|\mathfrak{F}^{-1}\big(\widehat{\mathcal{L}}^2_0(t,\xi)\big)(t,x)\big| &\lesssim |x|^{-m} \Big\|\mathfrak{F}^{-1}\Big(\partial_\xi^\alpha \big(\widehat{\mathcal{L}}^2_0(t,\xi) \big) \Big)(t,\cdot)\Big\|_{L^\ity} \\
&\lesssim |x|^{-m}\big\|\partial_\xi^\alpha \big(\mathcal{L}^2_0(t,\xi) \big)(t,\cdot) \big\|_{L^1} \\
&\lesssim |x|^{-m} e^{-ct} \int_1^{\ity} |\xi|^{-\sigma- m+n-1}d|\xi| \\
&\lesssim e^{-ct}\begin{cases}
|x|^{-(n-1)} &\text{ if } 0< |x| \le 1 \text{ and } m=n-1, \\
|x|^{-(n+1)} &\text{ if } |x| \ge 1 \text{ and } m=n+1,
\end{cases}
\end{align*}
where the assumption $\sigma>1$ comes into play. Hence, we arrive at
\begin{align*}
\big\|\mathfrak{F}^{-1}\big(\widehat{\mathcal{L}}^2_0(t,\xi)\big)(t,\cdot)\big\|_{L^1} &\lesssim \int_{|x|\le 1} \big|\mathfrak{F}^{-1}\big(\widehat{\mathcal{L}}^2_0(t,\xi)\big)(t,x)\big|dx \\
&\qquad+ \int_{|x|\ge 1} \big|\mathfrak{F}^{-1}\big(\widehat{\mathcal{L}}^2_0(t,\xi)\big)(t,x)\big|dx \\
&\lesssim e^{-ct}\int_0^1 d|x|+ e^{-ct}\int_1^\ity |x|^{-2}d|x| \lesssim e^{-ct}.
\end{align*}
Then, employing Young's convolution inequality we have proved the second statement in Proposition \ref{md2.21}. In the same way, we may also conclude the last statement and the third statement in Proposition \ref{md2.21}, respectively, by using (\ref{bd0.2.1.15}) and (\ref{bd0.2.1.17}). Let us come back to estimate the first statement. Indeed, we can see that
\begin{align}
\partial_t^{j}|D|^a \mathcal{K}^1_{u_0}(t,x) &= \partial_t^{j}|D|^a \mathfrak{F}^{-1}\Big(e^{\lambda_1(\xi)t}\chi_3(|\xi|) \widehat{u}_0(\xi) \Big)(t,x) \nonumber \\
&\quad + \partial_t^{j+1}|D|^a \mathfrak{F}^{-1}\Big(\frac{e^{\lambda_1(\xi)t}}{\lambda_2(\xi)- \lambda_1(\xi)}\chi_3(|\xi|) \widehat{u}_0(\xi) \Big)(t,x), \label{md0.2.2.1}
\end{align}
by using the relation
$$ \frac{\lambda_2(\xi) e^{\lambda_1(\xi) t}}{\lambda_2(\xi)- \lambda_1(\xi)}= e^{\lambda_1(\xi) t}+ \partial_t \Big(\frac{e^{\lambda_1(\xi) t}}{\lambda_2(\xi)- \lambda_1(\xi)}\Big). $$
In an analogous treatment to get the third statement, we derive the following estimate for the second term:
\begin{equation}
\Big\|\partial_t^{j+1}|D|^a \mathfrak{F}^{-1}\Big(\frac{e^{\lambda_1(\xi)t}}{\lambda_2(\xi)- \lambda_1(\xi)}\chi_3(|\xi|) \widehat{u}_0(\xi) \Big)(t,\cdot)\Big\|_{L^1} \lesssim e^{-ct} \|u_0\|_{H^{[a-\sigma]^+}_1}. \label{md0.2.2.2}
\end{equation}
In order to control the first term, using the relation $\lambda_1(\xi)= -1-\phi(\xi)$ we write
$$ e^{\lambda_1(\xi) t}= e^{-t} e^{-\phi(\xi)t}= e^{-t}- t e^{-t} \phi(\xi) \int_0^{^1} e^{-\phi(\xi) t r}dr. $$
Hence, we obtain
\begin{align}
&\mathfrak{F}^{-1}\Big(e^{\lambda_1(\xi)t}\chi_3(|\xi|) \widehat{u}_0(\xi) \Big)(t,x) \\
&\qquad= e^{-t}\mathfrak{F}^{-1}\big(\widehat{u}_0(\xi) \big)(x)- e^{-t}\mathfrak{F}^{-1}\big((1-\chi_3(|\xi|)) \widehat{u}_0(\xi)\big)(x) \nonumber \\
&\qquad \quad -t e^{-t}\mathfrak{F}^{-1}\Big(\phi(\xi) \chi_3(|\xi|) \widehat{u}_0(\xi) \int_0^{^1} e^{-\phi(\xi) t r}dr\Big)(t,x). \label{md0.2.2.3}
\end{align}
Obviously, we have
\begin{equation}
\Big\|\partial_t^{j}|D|^a \Big(e^{-t}\mathfrak{F}^{-1}\big(\widehat{u}_0(\xi)\big)\Big)(t,\cdot)\Big\|_{L^1}= e^{-t} \big\||D|^a u_0 \big\|_{L^1} \lesssim e^{-t}\|u_0\|_{H^a_1}. \label{md0.2.2.4}
\end{equation}
Now, we re-write
\begin{align*}
&\partial_t^{j}|D|^a \Big(t e^{-t}\mathfrak{F}^{-1}\Big(\phi(\xi) \chi_3(|\xi|) \widehat{u}_0(\xi)  \int_0^{^1} e^{-\phi(\xi) t r}dr\Big)\Big)(t,x) \\
&\quad = \sum_{\ell=0}^j \partial^{j-\ell}_t(t e^{-t})\,\partial_t^{\ell}|D|^a \mathfrak{F}^{-1}\Big(\phi(\xi) \chi_3(|\xi|) \widehat{u}_0(\xi)  \int_0^{^1} e^{-\phi(\xi) t r}dr\Big)(t,x) \\
&\quad= \sum_{\ell=0}^j \partial^{j-\ell}_t(t e^{-t})\,\mathfrak{F}^{-1}\Big(\phi^{\ell+1}(\xi)|\xi|^{\min\{a,\sigma\}} \chi_3(|\xi|) \int_0^{^1} e^{-\phi(\xi) t r}(-r)^\ell\,dr\Big)(t,x) \\
&\hspace{5cm} \ast  |D|^{[a-\sigma]^+}u_0(x) \\
&\quad =: \sum_{\ell=0}^j \partial^{j-\ell}_t(t e^{-t})\,\mathfrak{F}^{-1}\big(\widehat{\mathcal{L}}^1_0(t,\xi)\big)(t,x) \ast  |D|^{[a-\sigma]^+}u_0(x).
\end{align*}
Thanks to (\ref{bd0.2.1.11}) and (\ref{bd0.2.1.3}), by using the Leibniz rule we have
$$ \big|\partial_\xi^{\alpha}\big(\widehat{\mathcal{L}}^1_0(t,\xi)\big)\big| \lesssim e^{\frac{t}{2}}|\xi|^{-2\sigma \ell- 2\sigma+ \min\{a,\sigma\}- |\alpha|} \lesssim e^{\frac{t}{2}}|\xi|^{-\sigma -|\alpha|}. $$
Using again (\ref{Fouriertransform}), and carrying out $n-1$ and $n+1$ steps of partial integration imply
$$\big|\mathfrak{F}^{-1}\big(\widehat{\mathcal{L}}^1_0(t,\xi)\big)(t,x)\big| \lesssim e^{\frac{t}{2}}\begin{cases}
|x|^{-(n-1)} &\text{ if } 0< |x| \le 1, \\
|x|^{-(n+1)} &\text{ if } |x| \ge 1.
\end{cases} $$
It follows
$$\Big|\sum_{\ell=0}^j \partial^{j-\ell}_t(t e^{-t})\,\mathfrak{F}^{-1}\big(\widehat{\mathcal{L}}^1_0(t,\xi)\big)(t,x)\Big| \lesssim e^{-ct}\begin{cases}
|x|^{-(n-1)} &\text{ if } 0< |x| \le 1, \\
|x|^{-(n+1)} &\text{ if } |x| \ge 1,
\end{cases} $$
where $c$ is a suitable positive constant. Therefore, we derive
$$ \Big\|\sum_{\ell=0}^j \partial^{j-\ell}_t(t e^{-t})\,\mathfrak{F}^{-1}\big(\widehat{\mathcal{L}}^1_0(t,\xi)\big)(t,\cdot)\Big\|_{L^1} \lesssim e^{-ct}. $$
Applying Young's convolution inequality gives
\begin{align}
&\Big\|\partial_t^{j}|D|^a \Big(t e^{-t}\mathfrak{F}^{-1}\Big(\phi(\xi) \chi_3(|\xi|) \widehat{u_0}(\xi)  \int_0^{^1} e^{-\phi(\xi) t r}dr\Big)\Big)(t,\cdot)\Big\|_{L^1} \nonumber \\ 
&\qquad \lesssim e^{-ct}\|u_0\|_{H^{[a-\sigma]^+}_1}. \label{md0.2.2.5}
\end{align}
Moreover, due to $1- \chi_3 \in C_0^{\ity}$, we derive
$$ \Big\|\partial_t^{j}|D|^a \Big(e^{-t}\mathfrak{F}^{-1}\big(1- \chi_3(|\xi|)\big)\Big)(t,\cdot)\Big\|_{L^1} \lesssim e^{-t}. $$
By using again Young's convolution inequality we obtain
\begin{equation}
\Big\|\partial_t^{j}|D|^a \Big(e^{-t}\mathfrak{F}^{-1}\big((1- \chi_3(|\xi|))\widehat{u}_0(\xi)\big)\Big)(t,\cdot)\Big\|_{L^1}\lesssim e^{-t}\|u_0\|_{L^1}. \label{md0.2.2.6}
\end{equation}
Combining from (\ref{md0.2.2.1}) to (\ref{md0.2.2.6}) we may conclude the first statement in Proposition \ref{md2.21}. Summarizing, the proof of Proposition \ref{md2.21} is completed.
\end{proof}
\noindent From the statements in Proposition \ref{md2.21} we obtain immediately the following result.
\bmd
Let $\sigma>1$ and $n\ge 1$. The Sobolev solutions to (\ref{pt1.1}) satisfy the $L^1-L^1$ estimates
$$ \big\|\partial_t^{j}|D|^a u_{\chi_3}(t,\cdot)\big\|_{L^1} \lesssim e^{-ct}\Big(\|(u_0,u_1)\|_{H^{2\sigma j +[a- \sigma]^+}_1}+\|u_0\|_{H^a_1}+\|u_1\|_{H^{[a-\sigma]^+}_1}\Big), $$
where $c$ is a suitable positive constant, for any $t> 0$, $a \ge 0$ and for all integer number $j\ge 0$.
\label{md2.22}
\emd

\subsection{Estimates for middle frequencies}
Now let us turn to consider some estimates for middle frequencies, where $3^{-\frac{1}{\sigma}}\le |\xi| \le 3^{\frac{1}{\sigma}}$. Our goal is to clarify the exponential decay for solutions and some of their derivatives to (\ref{pt1.1}) localized to middle frequencies, which were neglected or not well-studied in the references.
\bmd
Let $\sigma>1$ and $n\ge 1$. The Sobolev solutions to (\ref{pt1.1}) satisfy the $L^1-L^1$ estimates
$$ \big\|\partial_t^{j}|D|^a u_{\chi_2}(t,\cdot)\big\|_{L^1} \lesssim e^{-ct}\|(u_0,u_1)\|_{L^1}, $$
where $c$ is a suitable positive constant, for any $t> 0$, $a \ge 0$ and $j= 0,1$.
\label{md2.3}
\emd
\begin{proof}
At first, with $3^{-\frac{1}{\sigma}}\le |\xi|\le 3^{\frac{1}{\sigma}}$ we use Cauchy's integral formula to re-write  the above Fourier multipliers in the following form:
\begin{align}
&\widehat{\mathcal{K}}_0(t,\xi)\chi_2(|\xi|)= \frac{1}{2\pi i}\Big(\int_{\Gamma}\frac{(z+|\xi|^{2\sigma})e^{zt}}{z^2+ |\xi|^{2\sigma}z+ |\xi|^{2\sigma}}dz\Big)\chi_2(|\xi|), \label{pro3.1.1.Midllezone} \\
&\widehat{\mathcal{K}}_1(t,\xi)\chi_2(|\xi|)= \frac{1}{2\pi i}\Big(\int_{\Gamma}\frac{e^{zt}}{z^2+ |\xi|^{2\sigma}z+ |\xi|^{2\sigma}}dz\Big)\chi_2(|\xi|), \label{pro3.1.2.Midllezone}
\end{align}
where $\Gamma$ is a closed curve containing the two characteristic roots $\lambda_{1,2}$. We can see that $\lambda_1= \lambda_2$ when $|\xi|= 2^{\frac{1}{\sigma}}$ and $\big\{\xi \in \R^n\,\, : \,\,|\xi|= 2^{\frac{1}{\sigma}}\big\}$ is not a singular set because we may give equivalent formulas as follows:
\begin{align*}
\widehat{\mathcal{K}}_0(t,\xi)= e^{\lambda_2 t}- \lambda_2 e^{\lambda_2 t}\int_0^t e^{(\lambda_1- \lambda_2) s}ds, \quad \widehat{K_1}(t,\xi)= e^{\lambda_2 t}\int_0^t e^{(\lambda_1- \lambda_2) s}ds.
\end{align*}
Therefore, it is reasonable to assume $\lambda_1 \neq \lambda_2$. Since $3^{-\frac{1}{\sigma}} < |\xi| < 3^{\frac{1}{\sigma}}$, this curve additionally is contained in $\{z\in \C\,\, : \,\, \text{\fontshape{n}\selectfont Re}\,z \le -c_0\}$, where $c_0$ is a positive constant.
In order to verify (\ref{pro3.1.1.Midllezone}) we express
$$ \frac{(z+|\xi|^{2\sigma})e^{zt}}{z^2+ |\xi|^{2\sigma}z+ |\xi|^{2\sigma}}= \frac{(z+|\xi|^{2\sigma})e^{zt}}{(z- \lambda_1)(z- \lambda_2)}= -\frac{\lambda_2}{\lambda_1- \lambda_2}\frac{e^{zt}}{z- \lambda_1}+ \frac{\lambda_1}{\lambda_1- \lambda_2}\frac{e^{zt}}{z- \lambda_2}. $$
For this reason, applying Cauchy's integral formula we obtain
\begin{align*}
&\frac{1}{2\pi i}\int_{\Gamma}\frac{(z+|\xi|^{2\sigma})e^{zt}}{z^2+ |\xi|^{2\sigma}z+ |\xi|^{2\sigma}}dz \\
&\qquad= -\frac{\lambda_2}{\lambda_1- \lambda_2}\Big(\frac{1}{2\pi i}\int_{\Gamma_1}\frac{e^{zt}}{z- \lambda_1}dz\Big)+ \frac{\lambda_1}{\lambda_1- \lambda_2}\Big(\frac{1}{2\pi i}\int_{\Gamma_2}\frac{e^{zt}}{z- \lambda_2}dz\Big) \\
&\qquad= -\frac{\lambda_2}{\lambda_1- \lambda_2}e^{\lambda_1 t}+ \frac{\lambda_1}{\lambda_1- \lambda_2}e^{\lambda_2 t}\chi_2(|\xi|)= \widehat{\mathcal{K}}_0(t,\xi)
\end{align*}
for middle frequencies. Here we split the curve $\Gamma$ into two closed sub-curves separately $\Gamma_1$ and $\Gamma_2$ containing $\lambda_1$ and $\lambda_2$, respectively. In the same way we may conclude the relation (\ref{pro3.1.2.Midllezone}). Now, taking account of estimates for $\widehat{\mathcal{K}}_0(t,\xi)$ we get
\begin{align*}
\mathfrak{F}^{-1}\big(|\xi|^a \widehat{\mathcal{K}}_0(t,\xi) \chi_2(|\xi|)\big)&= \int_{\R^n}e^{ix\xi}|\xi|^a \widehat{\mathcal{K}}_0(t,\xi) \chi_2(|\xi|)d\xi \\
&= \sum_{k=1}^n \frac{x_k}{i|x|^2}\int_{\R^n} \partial_{\xi_k}\big(e^{ix\xi}\big) |\xi|^a \widehat{\mathcal{K}}_0(t,\xi) \chi_2(|\xi|) d\xi,
\end{align*}
where we used (\ref{Fouriertransform}). By induction argument, carrying out $m$ steps of partial integration we derive
$$ \mathfrak{F}^{-1}\big(|\xi|^a \widehat{\mathcal{K}}_0(t,\xi) \chi_2(|\xi|)\big)= C \sum_{|\alpha|=m} \frac{(ix)^\alpha}{|x|^{2|\alpha|}} \mathfrak{F}^{-1}\Big(\partial_\xi^\alpha \big(|\xi|^a \widehat{\mathcal{K}}_0(t,\xi)\chi_2(|\xi|)\big) \Big), $$
for any non-negative integer $m$ and $C$ is a suitable constant. Hence, we arrive at the following estimates:
\begin{align*}
\big|\mathfrak{F}^{-1}\big(|\xi|^a \widehat{\mathcal{K}}_0(t,\xi)\chi_2(|\xi|)\big)\big| &\lesssim |x|^{-m} \Big\|\mathfrak{F}^{-1}\Big(\partial_\xi^\alpha \big(|\xi|^a \widehat{\mathcal{K}}_0(t,\xi)\chi_2(|\xi|)\big) \Big)\Big\|_{L^\ity}\\
&\lesssim |x|^{-m}\big\|\partial_\xi^\alpha \big(|\xi|^a \widehat{\mathcal{K}}_0(t,\xi)\chi_2(|\xi|)\big) \big\|_{L^1} \lesssim |x|^{-m} e^{-ct},
\end{align*}
where $c$ is a suitable positive constant, since $3^{-\frac{1}{\sigma}}< |\xi|< 3^{\frac{1}{\sigma}}$. This estimate immediately implies
$$ \big\|\mathfrak{F}^{-1}\big(|\xi|^a \widehat{\mathcal{K}}_0\chi_2(|\xi|)\big)(t,\cdot)\big\|_{L^1} \lesssim e^{-ct}. $$
In an analogous way we may also conclude
$$ \big\|\mathfrak{F}^{-1}\big(|\xi|^a \widehat{\mathcal{K}}_1\chi_2(|\xi|)\big)(t,\cdot)\big\|_{L^1} \lesssim e^{-ct}. $$
Similarly, we may arrive at the exponential decay for the following estimates:
\begin{align*}
&\big\|\mathfrak{F}^{-1}\big(|\xi|^a \partial_t \widehat{\mathcal{K}}_0\chi_2(|\xi|)\big)(t,\cdot)\big\|_{L^1}\lesssim e^{-ct}, \\
&\big\|\mathfrak{F}^{-1}\big(|\xi|^a \partial_t \widehat{\mathcal{K}}_1\chi_2(|\xi|)\big)(t,\cdot)\big\|_{L^1}\lesssim e^{-ct},
\end{align*}
where we notice that
$$ \partial_t \widehat{\mathcal{K}}_0(t,\xi)= -|\xi|^{2\sigma}\widehat{\mathcal{K}}_1(t,\xi) \quad \text{ and }\quad \partial_t \widehat{\mathcal{K}}_1(t,\xi)= \widehat{\mathcal{K}}_0(t,\xi)- |\xi|^{2\sigma}\widehat{\mathcal{K}}_1(t,\xi). $$
Therefore, applying Young's convolution inequality we get
\begin{align*}
\big\|\partial_t^{j}|D|^a u_{\chi_2}(t,\cdot)\big\|_{L^1} &\lesssim \big\|\mathfrak{F}^{-1}\big(|\xi|^a \partial^j_t \widehat{\mathcal{K}}_0\chi_2(|\xi|)\big)(t,\cdot)\big\|_{L^1} \|u_0\|_{L^1} \\
&\qquad+ \big\|\mathfrak{F}^{-1}\big(|\xi|^a \partial^j_t \widehat{\mathcal{K}}_1\chi_2(|\xi|)\big)(t,\cdot)\big\|_{L^1} \|u_1\|_{L^1} \\
&\lesssim e^{-ct} \|(u_0,u_1)\|_{L^1}.
\end{align*}
Summarizing, the proof to Proposition \ref{md2.3} is completed.
\end{proof}
\begin{proof}[Proof of Theorem \ref{dl1.1}]
We combine the statements from Propositions \ref{md2.1}, \ref{md2.22} and \ref{md2.3} to conclude immediately all the desired estimates. This completes our proof.
\end{proof}


\end{document}